\documentclass[11pt]{article}
\usepackage{latexsym}
\usepackage{graphicx}
\usepackage{epstopdf}
\usepackage{lineno}
\usepackage[sort]{natbib}
\newcommand{\coln}{\hspace*{-6pt}{\bf :}}
\newtheorem{theorem}{Theorem}
\newtheorem{lemma}{Lemma}
\newtheorem{corol}{Corollary}
\newtheorem{property}{Property}

\newcommand{\Theorem}[1]{\begin{theorem}\coln ~#1 \end{theorem}}
\newcommand{\Lemma}[1]{\begin{lemma}\coln ~#1 \end{lemma}}

\newcommand{\proof}[1]{\noindent{\bf Proof:~}#1\hfill ~$\Box$}
\newcommand{\be}{\begin{equation}}
\newcommand{\ee}{\end{equation}}

\newcommand{\narrow}{\setlength{\itemsep}{1pt}\setlength{\parskip}{6pt}
\setlength{\parsep}{0pt}
\renewcommand{\baselinestretch}{1}
\large
\normalsize}
\usepackage{color}

\begin{document}


\title{\bf Directional Approach to Gradual Cover: The Continuous Case}
\author{Tammy Drezner, Zvi Drezner and Pawel Kalczynski\\
Steven G. Mihaylo College of Business and Economics\\
California State University-Fullerton\\
Fullerton, CA 92834.\\e-mail: tdrezner@fullerton.edu; zdrezner@fullerton.edu; pkalczynski@fullerton.edu}
\date{}
\maketitle

\begin{abstract}
The objective of the cover location models is covering demand by facilities within a given  distance. The gradual (or partial) cover replaces abrupt drop from full cover to no cover by defining gradual decline in cover.  In this paper we use a recently proposed rule for calculating the joint cover of a demand point by several facilities termed ``directional gradual cover".
Contrary to all gradual cover models, the joint cover depends on the facilities' directions. In order to calculate the joint cover, existing models apply the partial cover by each facility disregarding their direction. We develop a genetic algorithm to solve the facilities location problem and also solve the problem for facilities that can be located anywhere in the plane.
The proposed modifications were extensively tested on a case study of covering Orange County, California.
\end{abstract}


\renewcommand{\baselinestretch}{1.6}
\renewcommand{\arraystretch}{0.625}
\large
\normalsize

\section{Introduction}

Cover location models constitute a main branch of location analysis. A demand point is covered by a facility within a certain distance \citep{CR74,RTF76}. A given number of facilities need to be located in an area so as to provide as much cover as possible. Such models are used for modeling cover provided by emergency facilities such as ambulances, police cars, or fire trucks. They are also used to model cover by transmission towers for cell-phone, TV, radio, radar among others.

In gradual cover models (also referred to as partial cover) it is assumed that up to a certain distance $r$ the demand point is fully covered and beyond a greater distance $R$ it is not covered at all. Between these two extreme distances the demand point is partially covered. There exist several formulations of the partial cover. \citet{BK02} suggested a declining step function between $r$ and $R$. \citet{DWD04} suggested a linear decline in cover between $r$ and $R$, and \citet{DDG09} suggested a linear decline between random values of $r$ and $R$. \citet{DDK17} proposed that the demand point is represented by a circle of diameter $R-r$ and the cover is the intersection area of the circle centered at the demand point and the circle centered at the facility with radius $\frac{R+r}2$. The original cover models are a special case of gradual cover models. When $R=r$, cover drops abruptly from full cover to no cover. All the applications listed above for standard cover models are modeled more realistically by gradual cover. In reality, cover does not drop abruptly at a certain distance.

\citet{CR84} were the first to propose a discrete gradual cover model.
The network version with a step-wise cover function is discussed in \cite{BK02}. The network and discrete models with a general non-increasing cover function were analyzed in \cite{BKD03}. The  single-facility planar model  with a linearly decreasing cover function between the distance of full coverage and the distance of no coverage was optimally solved  in \cite{DWD04}. The stochastic version was analyzed and optimally solved in \cite{DDG09}.  Additional references include \cite{KK04,EM09,DrDr14,BDK11b,BDK10}. For a review of cover models see \cite{Plas02,GM15,S11,DDK17,BDK10}.

A main issue when $k$ facilities partially cover a demand point is estimating the total cover, also referred to as joint cover.  One approach is to interpret partial cover  as the probability $p_j$ of cover for $j=1,\ldots,k$. If the probabilities are independent of each other, the joint cover is $1-\prod\limits_{j=1}^k(1-p_j)$.
For a discussion of ways to estimate the joint cover see \cite{DrDr14,BDK11b,EM09,KK04}. Existing models just apply the partial cover by each facility disregarding their direction. Two facilities located north of the demand point (which is actually an area and not a point) at distances of 1 and 2 miles provide the same joint cover as two facilities located one mile north and 2 miles south of the demand point. However, when both facilities are north of the demand point there is considerable overlap in cover while if one is to the north and one to the south, the possible overlap is smaller. Therefore, the joint cover when both facilities are at the same direction is lower than the cover if they are located in opposite directions.

To address this issue of direction, \citet{DDK17} introduced a directional approach to gradual cover. Customers at demand points reside in areas represented by circular discs rather than mathematical points and the facility covers points within a given distance. To find a demand point's partial cover by one facility, the intersection area of two discs (the demand point's disc and the facility's coverage disc) is calculated. The partial cover of the demand point is the intersection area divided by the demand point's disc area.

\begin{figure}[ht!]
	\begin{center}
		\setlength{\unitlength}{1in}
		\begin{picture}(3,2.8)
		\includegraphics[width=3in]{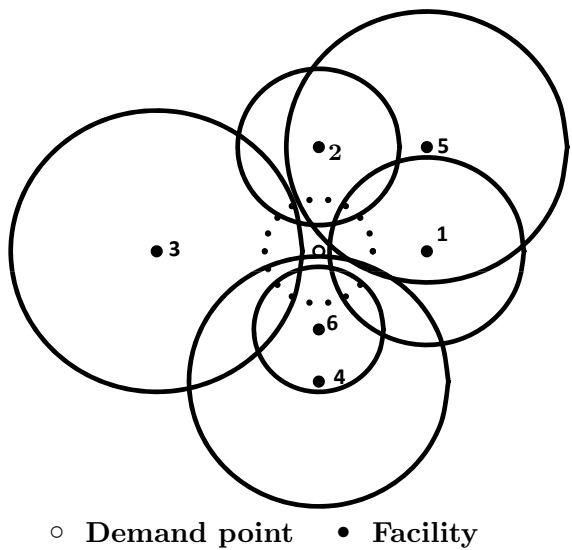}
		\put(-2.7,-0.1){\circle{0.05}}
		\put(-2.55,-0.15){\bf Demand point}
		\put(-1.2,-0.1){\circle*{0.05}}
		\put(-1.05,-0.15){\bf Facility}
		\put(-1.28,1.85){\footnotesize\bf 2}
		\end{picture}
		\bigskip
		\caption{\label{example}An Illustration of Directional Cover}
	\end{center}
\end{figure}

If several facilities exist in an area, the joint cover of a demand point is the union of the individual areas covered by the facilities. This joint cover depends on the distances to the facilities and on their directions. See Figure \ref{example} for an example of joint cover. Six facilities are located in the area. Their coordinates and cover radii $D_j$ are depicted in Table \ref{fac}. The demand point is located at (0, 0). The disc defined by the ``demand point" of radius 1 is marked by dots. The discs centered at the facilities are marked with thick circular circumferences.

\begin{table}[ht!]
	\begin{center}										
		\caption{\label{fac}Parameters of Six Facilities}
		\medskip															\begin{tabular}{|c||c|c|c|}															
			\hline
			$j$&$x_j$&$y_j$&$D_j$\\
			\hline
			1&	2&	0&	1.8\\
			2&	0&	2&	1.5\\
			3&	-3&	0&	2.7\\
			4&	0&	-2.5&	2.4\\
			5&	2&	2&	2.6\\
			6&	0&	-1.5&	1.2\\
			\hline
		\end{tabular}					
	\end{center}
\end{table}

Each of the six facilities covers part of the demand point's disc. The joint cover of the demand point is the union of the intersection areas between the facilities' discs and the disc surrounding the demand point.

The joint cover of several facilities depends on the direction of the facilities. To illustrate the effect of direction, consider the  example in Figure \ref{example} of cover by two facilities excluding the other four facilities in the area. Facility 2 is to the north of the demand point and facilities 4 and 6 are south of the demand point. The total area covered by facilities 2 and 4 is the sum of the cover areas. Facility 2 covers part of the northern region of the neighborhood while facility 4 covers part of the southern region. The covered areas do not intersect. The same is true for facilities 2 and 6 and also for facilities 1 and 3. However, facilities 4 and 6 are both south of the demand point and cover the southern part of the neighborhood. Because of the overlap of the cover areas the total cover is the area covered by facility 4 while no additional cover is provided by facility 6. The total cover of facilities 1 and 4 is less than the sum of the areas but more than the larger area. There is an intersection area of the two cover areas which is counted only once. Contrary to all gradual cover models, the joint cover depends on the direction. In order to calculate the joint cover, it is not sufficient to only have the value of its partial cover for each facility.

The main contribution of the paper is solving the model when facilities can be located anywhere in the plane and not on a given set of locations. As can be expected, such a solution provides a better cover. In order to evaluate the value of the objective function anywhere in the plane for standard non-linear non-convex solvers (we applied SNOPT), we applied the hexagonal pattern numerical integration approach also used in \cite{DDS18b}.
We also solve the discrete directional cover model by a genetic algorithm getting better results than those obtained in \cite{DDK17} by other heuristic approaches (Ascent, Tabu search, and Simulated annealing).

The paper is organized as follows. In the next section we detail two approaches to estimate joint cover by numerical integration. In Section \ref{gen} we propose a genetic algorithm for solving the discrete problem where a set of potential locations for the facilities is given. In Section \ref{sec4} we solve the continuous problem where the facilities are located anywhere in the plane. Two solution approaches are suggested: (i) applying general purpose solvers and (ii) a Nelder Mead \citep{NM65,DW87} approach. Extensive computational experiments are performed on a case study of covering Orange County, California and  presented in Section \ref{case}. We conclude the paper in Section \ref{concl}.

\section{\label{calc}Calculating the Joint Cover Numerically}

Evaluating the union of the individual cover areas is complex. The area can be calculated by two-dimensional integration. Two numerical integration methods for approximating the cover area are used. Gaussian quadrature is used in \citet{DDK17} and the calculations details are available there. The hexagonal pattern is also used in \citet{DDS18b} and details are summarized in the appendix.

\begin{description}
	\item[Gaussian Quadrature:] 
The percent of coverage of a circle can be accurately calculated by finding the covered arcs on the circumference of that circle.  \citet{DDK17} used this observation to estimate the union's area by Gaussian quadrature based on Legendre polynomials \citep{AS72}. The  {\it exact} integral over a circle is found and each cover by a circle is multiplied by the appropriate weight given in Table \ref{legend}.

\item[Hexagonal Pattern:] The union's area is estimated by evaluating the cover at points in an hexagonal pattern. This hexagonal pattern is easily generated as detailed in \citet{DDS18b} and in the appendix.

\end{description}

\begin{table}[ht!]
	\begin{center}															
		\caption{\label{legend}Adjusted Legendre-Gaussian Quadrature Parameters}
		\medskip															\begin{tabular}{|c||c|c|}
			\hline
			$j$&$u_j$&$w_j$\\
			\hline
			1&	0.1142223084&	0.0333356722\\
			2&	0.2597466394&	0.0747256746\\
			3&	0.4003688498&	0.1095431813\\
			4&	0.5322614986&	0.1346333597\\
			5&	0.6523517690&	0.1477621124\\
			6&	0.7579163341&	0.1477621124\\
			7&	0.8465800004&	0.1346333597\\
			8&	0.9163540714&	0.1095431813\\
			9&	0.9656768007&	0.0747256746\\
			10&	0.9934552150&	0.0333356722\\
			\hline
		\end{tabular}
	\end{center}
\end{table}

\begin{figure}[ht!]
	\begin{center}
		\setlength{\unitlength}{1in}
		
		\begin{picture}(6.5,3.1)
		\thicklines
		\includegraphics[width=3in]{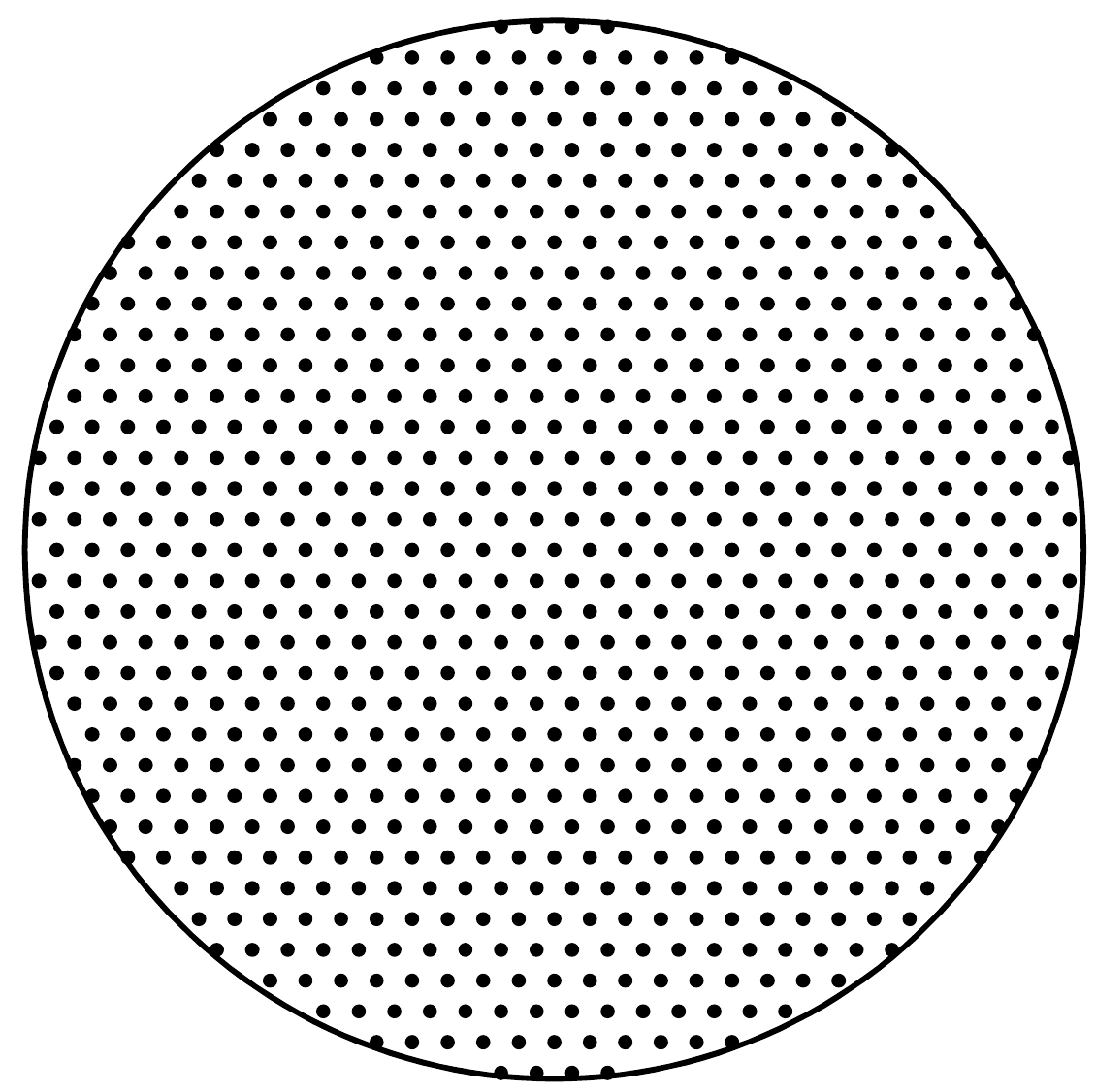}
		\includegraphics[width=3in]{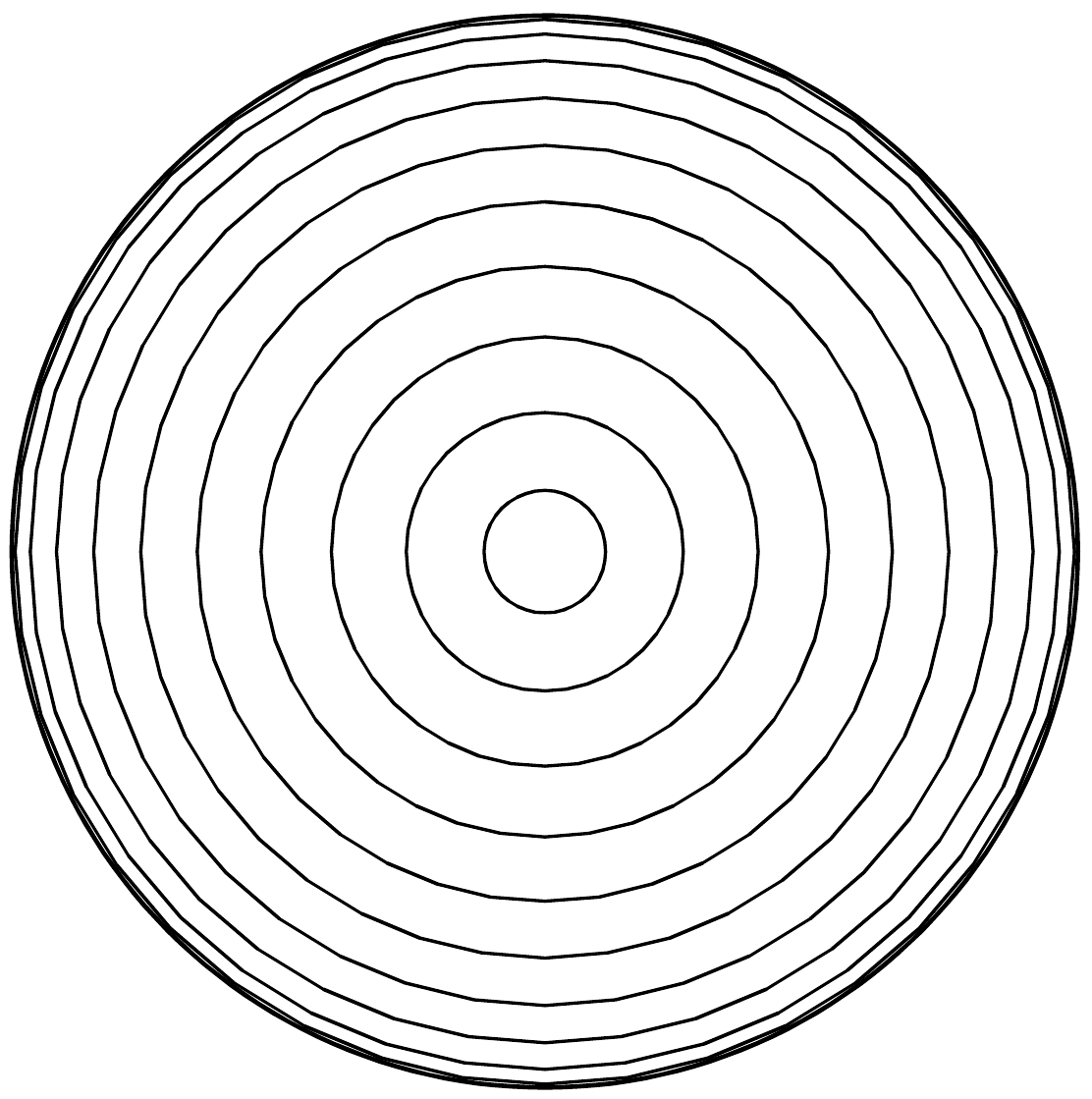}
		\put(-5.8,-0.25){\bf Hexagonal Pattern, 805 Points}
		\put(-3.1,-0.25){\bf Gaussian Quadrature, 10 Concentric Circles}
		\put(-4.527,1.483){\circle{0.1}}
			\put(-4.527,1.483){\circle{0.08}}
		\put(-1.503,1.5){\circle{0.1}}
				\put(-1.503,1.5){\circle{0.08}}
		
		\put(-3.5,-0.4){\circle{0.08}}
		\put(-3.5,-0.4){\circle{0.07}}
			\put(-3.4,-0.45){{\bf Circle's Center}}
		\end{picture}
		\medskip		
		\caption{\label{hex}Comparing Numerical Integration Methods}
	\end{center}
\end{figure}

The hexagonal pattern for $N=805$ points is depicted on the left panel of Figure~\ref{hex}. The largest distance between the points and the center of the circle is 0.99342. The ten points (which are circles) Gaussian quadrature are depicted on the right panel of Figure \ref{hex}. The largest circle has a radius of 0.99346 which is not distinguishable to the naked eye from the circle of radius 1 surrounding the demand point.

We calculated the total cover by the demand point in an example problem (Table~\ref{fac}, Figure~\ref{example}) with demand point's radius between 1.0 and 2.0 by a simulation of a billion points randomly generated in the circle. Each instance required less than a minute of computer time. The standard error of the results is about $1.4\times10^{-5}$. The results by the simulation, the Gaussian quadrature, and hexagonal patterns for various number of points, are depicted in Table \ref{area}.

\begin{table}[ht!]
	\begin{center}										
		\caption{\label{area}Cover Estimations by Demand Point Radius}
		\medskip															\begin{tabular}{|c||c||c|c|c|c|}
			\hline
			$R$&Sim.&Gauss&$N=199$&$N=397$&$N=805$\\													
			\hline
			1.0&	0.920&	0.923&	0.925&	0.917&	0.924\\
			1.1&	0.934&	0.933&	0.935&	0.932&	0.934\\
			1.2&	0.945&	0.947&	0.950&	0.937&	0.945\\
			1.3&	0.953&	0.954&	0.955&	0.960&	0.954\\
			1.4&	0.959&	0.960&	0.955&	0.962&	0.958\\
			1.5&	0.965&	0.965&	0.965&	0.965&	0.964\\
			1.6&	0.969&	0.968&	0.975&	0.970&	0.968\\
			1.7&	0.972&	0.970&	0.975&	0.975&	0.969\\
			1.8&	0.975&	0.973&	0.980&	0.975&	0.976\\
			1.9&	0.978&	0.976&	0.980&	0.977&	0.980\\
			2.0&	0.980&	0.978&	0.980&	0.977&	0.981\\
			\hline
			\multicolumn{2}{|l||}{Average$^*$}&0.0016&	0.0029&	0.0028&	0.0016\\
			\hline
			\multicolumn{6}{l}{$^*$ Average absolute deviation from the simulation.}		\end{tabular}					
	\end{center}
\end{table}

It should be noted that the results in Table \ref{area} are for only  one demand point. In our computational experiments we implemented a case study of 577 demand points. Therefore, the estimation of the total area covered is much more accurate (its standard error is about  $\frac{1}{\sqrt{577}}$ of the standard deviation of an individual demand point.) Moreover, when the number of facilities is small, many of the demand points are not covered at all and their estimation is also zero which is accurate. When $p$ is large, many demand points are fully covered and their estimation is also accurate with a value of 1.

\section{\label{gen} The Proposed Genetic Algorithm for the Discrete Problem}

\citet{DDK17} proposed and tested the ascent algorithm, tabu search \citep{GL97}, and simulated annealing \citep{KGV83} for the solution of this problem. Simulated annealing performed best. We construct and test the following genetic algorithm that performed much better than simulated annealing and tabu search. We experimented with many variants and parameters of the genetic algorithm and the following one performed the best in our experiments. We borrowed ideas used in \cite{ADE03,DM03,Dr03}. The value of the objective function can be calculated by the two methods detailed in Section \ref{calc}. We opted to apply the Gaussian quadrature procedure that was used in \cite{DDK17}.

Consider a problem with $n$ demand points and $m$ potential facilities' locations. Co-location is not beneficial for cover. When two facilities are located at the same point, the intersection area for the facility with a smaller or equal cover radius is included in the intersection area of the facility with the larger radius. The union is the intersection area of the facility with a larger cover radius. Therefore, the objective is to select the best $p$ of the $m$ potential locations for locating facilities without repetition.  Once the cover by selecting $p$ potential locations can be calculated, the problem reduces to selecting the best set of $p$ out of $m$ potential locations. If $p$ and $m$ are relatively small, total enumeration or a branch and bound algorithm can be used. Otherwise, a heuristic approach is necessary.

A list of selected locations is given. We define:
\begin{description}
\item[An Ascent Algorithm:]	Each iteration the value of the objective function is evaluated by replacing a selected location with a non-selected location. The best exchange is performed if it is an improving one. The procedure is stopped with a final solution when no improving solution is found.
\item[Reverse Greedy Process:] There are more than $p$ selected locations.  Each iteration the location whose removal decreases the objective function the least is removed until the number of selected locations is reduced to $p$.
\item[Restricted Ascent:] An ascent algorithm where the list of non-selected locations is restricted to a given set and does not necessarily include all non-selected locations. 
\end{description}

\subsection*{Outline of the Genetic Algorithm}

Every population member is a list of $p$ locations where facilities are located. The algorithm is run for a pre-specified number of generations, $g$.
\begin{enumerate} \narrow
\item A population of size $pop$ (we used $pop=100$) is randomly generated. Set $iter=0$.
\item \label{gst2}20\% of the population members are randomly selected and are possibly improved by an ascent algorithm.
\item\label{gst4} Two population members are selected by the parents selection rule detailed below.
\item An offspring is produced by the merging process detailed below.
\item If the offspring is worse than the worst population member or is identical to an existing population member, go to Step \ref{gst7}.
\item The offspring replaces the worst population member.
\item\label{gst7}Set $iter=iter+1$. If $iter\le g$ go to Step \ref{gst4}.	
\item The best population member is the result of the genetic algorithm.
\end{enumerate}

The following parents selection rule was proposed in \cite{DM03}. It is based on the biological concepts of inbreeding and outbreeding depression \citep{E07, FG00}. A successful offspring is formed when the parents are not too close genetically (inbreeding depression) or too dissimilar (outbreeding depression).

\subsubsection*{The Parents Selection Rule}
\noindent Define:
\begin{itemize} \narrow
\item[$P$] Number of potential second parents (a parameter).
\item[$c$] Similarity count. The number of common locations between two specific population members.
\end{itemize}
\noindent The process:
\begin{enumerate} \narrow
	\item One population member is randomly selected as the first parent.
	\item $P\ge 1$ potential second parents are randomly selected from the remaining $pop-1$ population members.
	\item For each potential second parent, the similarity count with the first parent, $c$, is found.
	\item The potential second parent with the smallest value of $c$ is selected.
\end{enumerate}

 An important part of a genetic algorithm is the merging process of two parents producing an offspring. The following merging process was found to be superior to many others that were tested.
 
\subsubsection*{The Merging Process}
\begin{enumerate} \narrow
\item \label{mst1}Two lists of locations are formed. The first list consists of all locations common to the two parents. The second list consists of all locations which exist only in one parent. The first list has $0\le c\le p$ members and the second list has $2p-2c$ members.
\item $p-c$ members are randomly selected from the second list and added to the first list creating a solution with $p$ locations. There are  $p-c$ remaining locations in the second list.
\item\label{step3} An ascent algorithm is applied on the created solution considering only exchanges between locations in the solution with locations in the second list. Each iteration of this ascent algorithm requires $p(p-c)$ evaluations of the objective function. If the solution is improved, the solution and the second list are updated and another iteration is performed until no improvement is found.
\item \label{step4}Once the result of the ascent algorithm is obtained, $\lfloor {\frac p2} 
\rfloor$ potential locations are randomly selected from the remaining $m-p$ potential locations creating a new second list.
\item \label{step4a}The ascent algorithm is repeated with the newly formed second list yielding the offspring.
\end{enumerate}

Note that as the number of potential second parents, $P$, increases, the similarity measure $c$ tends to decrease. Therefore, the number of objective function evaluations, $p(p-c)$, in Step \ref{step3} of the merging process is expected to increase leading to some increase in run times.

\subsection{Variants Tested}
The following describes successive improvements in the genetic algorithm that were tested.
\begin{enumerate} \narrow
\item Step \ref{gst2} of the genetic algorithm was added later.
\item We started with the merging process suggested in \cite{ADE03}. The two lists are created as in Step \ref{mst1} of the merging process. The two lists are combined into one. The locations in the first list stay in the solution and locations from the second list are removed from the combined list by a reverse greedy process.
\item We then tried performing an ascent algorithm on each offspring but this took too much computer time and necessitated a significant reduction in the number of generations.
\item  We then tried Step \ref{gst2} of the genetic algorithm on 100\% of the initial population members but it did not do well because there were many identical population members.
\item  We then changed the 100\% to 20\% as is described in Step \ref{gst2} of the genetic algorithm.
\item  We then tried the restricted ascent as detailed in Step \ref{step3} of the merging process getting relatively good results.
\item Steps \ref{step4}-\ref{step4a} were added to the merging process to prevent a situation where some locations are not present in any population member. Without these steps, such locations will not be present in any offspring. Consequently, they will not be in the final solution even if some of them are present in a better one. These two steps can be viewed as creating a mutation of some offspring \citep{SP94,H15,FSH08}.
\end{enumerate}

\section{\label{sec4} Finding a Solution Anywhere on the Plane}
 We first show that the facilities in an optimal solution anywhere in the plane are located within the convex hull of the demand points.

\Lemma{\label{conv}
An optimal location of one facility exists in the convex hull of the demand points.}
\proof{Let $Y$ be a location outside the convex hull. We show that there is a location $Z$ in the convex hull such that the cover by a facility located at $Z$ is at least the cover by a facility located at $Y$. \citet{WH73} showed that there is a point $Z$ in the convex hull (or any convex set) which is closer than point $Y$ to every point $X$ in the convex hull. Therefore, the partial cover of point $X$ by a facility located at $Z$ is at least the cover by a facility located at $Y$.}
 
\Theorem{\label{Th2}There exists an optimal solution to the multiple facility location problem where all facilities are located in the convex hull of the demand points.}
\proof{Similar to the arguments in \cite{HPT80}, suppose that a facility $Y$ is located outside the convex hull. By Lemma \ref{conv} there is a location $Z$ in the convex hull for that facility, while holding the other facilities in their locations, which cover at least the cover by a facility located at $Y$. The theorem follows by mathematical induction.
}
 
\subsection{Heuristic Algorithms}

Heuristic algorithms are designed for the location of facilities anywhere in the plane. The value of the objective function for any set of $p$ facilities can be calculated by numerical integration. Once a starting solution is selected, the solution is improved either by applying all purpose non-linear solvers or by a specially designed Nelder-Mead algorithm detailed below. Note that by Theorem \ref{Th2} there is no need for constraints.

We suggest two approaches for generating starting solutions. One approach is to randomly select $p$ demand points, and another is 
to find a solution to the discrete cover problem and use it as a starting solution. Finding the discrete solution requires extensive run time and does not provide a fair comparison with random starting solutions. Results with the discrete starting solution are reported because in many cases they provide the best known solution.

We propose to move the facilities to nearby locations that provide greater total cover. A facility located at a demand point completely covers that demand point. However this particular demand point is still fully covered if the facility is moved up to a distance $D-R$ from that demand point. The gradient of the coverage of the particular demand point is zero within the disc of radius $D-R$. However, a gradient search in the neighborhood of that demand point may increase total coverage of other demand points.

There is no analytical expression for the total cover and thus applying non-linear solvers for finding a good solution may be complex unless an hexagonal pattern is used. Multipurpose non-linear solvers are applied by providing a code for the evaluation of the objective function. For the Nelder-Mead algorithm we propose,
similarly to \citet{Co63,Co64}, to iteratively find a better location for one facility at a time (in random order) while holding the other $p-1$ facilities in their place, until convergence. Demand points that are fully covered by the $p-1$ fixed facilities can be removed from the optimization problem of locating a single facility.
Global optimization techniques such as ``Big Square Small Square" (BSSS) \citep{HPT81} or ``Big Triangle Small Triangle" (BTST) \citep{DS04} can be used for solving the single facility location problem. Since the algorithm is a heuristic, it is unnecessary to find the location of one facility optimally. We used the
 Nelder-Mead \citep{NM65,DW87} for relocating one facility because it is faster than optimal approaches.

\subsection{Nelder-Mead (for Maximization)}
Suppose that there are $k$ variables in a function to be maximized. In the original Nelder-Mead algorithm, a starting set of $K=k+1$ solutions forming a ``simplex" is generated.
We solve relocation of one facility while holding the other $p-1$ facilities fixed in their locations. Therefore, we have $k=2$ variables leading to a simplex of $K=3$ vertices. We select the present location as one of the vertices and randomly generate two additional vertices in a square centered at the present location.
We also experimented with a simplex of $K=4$ solutions. The simplex is moved by replacing the worst solution by a better one.

The Nelder-Mead algorithm is based on three parameters $\alpha,\beta$ and $\gamma$.
The recommended values are: for the reflection $\alpha=1$, for the contraction: $\beta=0.5$, and for the expansion $\gamma=2$. The reduction parameter that was 0.5 in the original Nelder-Mead algorithm can also be modified.

\subsubsection{Definitions}
\begin{itemize}\narrow
\item[$P_i$] for $i=1,\ldots,K$ are the solution points (the vertices of the simplex).
\item[$ÙF_i$] for $i=1,\ldots,K$ are the corresponding values of the objective function.
\item[$P_{\ell}$] is the worst solution point ($\ell$ for low).
\item[$P_{h}$] is the best solution point ($h$ for high).
\item[$\overline P$] is the average of all solution points excluding $P_{\ell}$.
\item[$F_s$] is the value of the objective function of the second highest solution point.
\end{itemize}

\subsubsection{The Nelder-Mead Algorithm}

\begin{enumerate} \narrow
	\item \label{st1} Determine the solutions at the vertices $P_{\ell},P_{h},\overline P$.
	\item Determine the values $F_{\ell},F_s,F_h$.
	\item If $F_h-F_{\ell}<\epsilon$, stop with $P_{h}$ as the solution.
	\item \label{st4}{\it Reflection}: Calculate $P^r=(1+\alpha)\overline P -\alpha P_{\ell}$ and its objective function value $F^r$.
	\item \label{st5}If $F_s\le F^r \le F_h$, replace $P_{\ell}$ by $P^r$  and return to Step \ref{st1}.
	\item \label{st6}If $F^r < F_{s}$ then define $P^t$: If $F^r>F_{\ell}$, $P^t=P^r$; else $P^t=P_{\ell}$.
	\item \label{st6a}{\it Contraction}: Calculate $P^c=\beta P^t+(1-\beta)\overline P$ and its objective value $F^c$.
	\begin{enumerate}
		\item If $F^c\ge F_s$, replace $P_{\ell}$ by $P^c$ and go to Step \ref{st1}.
		\item {\it Reduction}: If $F^c< F_s$, replace all points $P_i=(P_i+P_{h})/2$ and go to Step \ref{st1}.
	\end{enumerate}
\item \label{st7}{\it Expansion}: If $F^r>F_{h}$, calculate $P^e=(1+\gamma)\overline P -\gamma P_{\ell}$ and its objective $F^e$.
		\begin{enumerate}
			\item If $F^e\ge F_{h}$, replace $P_{\ell}$ by $P^e$ and go to Step \ref{st1}.
			\item Else, replace $P_{\ell}$ by $P^r$ and  go to Step \ref{st1}.
		\end{enumerate}
\end{enumerate}

\section{Case Study: Transmission Towers in Orange County, California} \label{case}

We investigate covering Orange County, California with transmission towers of cell phone, TV or radio.
The data from the 2000 census for Orange County, California is given in \cite{Dr04} and was also used in \cite{DrDr07,BDK09,DrDr14,BDK11b,DDK17,DDS06}. There are  577 census tracts and their population counts are given. The total Orange County population is 2,846,289.

Experiments using  SNOPT \citep{SNOPT} were run on a virtual server with 16 vCPUs and 128 GB of vRAM implemented in Matlab R2016b. An analytical gradient for the objective function is not available. The derivatives of the objective function are discontinuous. Therefore, the SNOPT non-linear all purpose solver is not effective.

For the Nelder-Mead approach, computer programs were coded in Fortran using double precision arithmetic and were compiled by an Intel 11.1 Fortran Compiler  using one thread with no parallel processing. They were run on a desktop with the Intel i7-6700 3.4GHz CPU processor and 16GB RAM.

Full cover within 2 miles and no cover beyond 4 miles were applied in \cite{BDK11b,DrDr14,DDK17}. To have comparable results we assign a radius of $R_i=1$ mile for each demand point  and a cover radius of $D_j=3$ miles for each tower. For each community with population $w_i$, the  proportion cover $0\le p_i\le 1$ is estimated by numerical integration. The proportion cover of all $n$ communities $0\le  \hat p \le 1$ is
$$
\hat p =\frac{\sum\limits_{i=1}^n w_ip_i}{\sum\limits_{i=1}^nw_i}
$$

\begin{table}[ht!]
	\begin{center}															
		\caption{\label{comp1}Comparing the Genetic Algorithm with Simulated Annealing ($p\le 22$)}
		\medskip															\begin{tabular}{|c|c||c|c||c|c||c|c||c|c|}															\hline
			&Best&\multicolumn{2}{|c||}{Gen ($P=1$)}&\multicolumn{2}{|c||}{Gen ($P=2$)}&\multicolumn{2}{|c||}{Gen ($P=3$)}&\multicolumn{2}{|c|}{S.A.}\\
			\cline{3-10}
			$p$&Known&$\dagger$&$\ddagger$&$\dagger$&$\ddagger$&$\dagger$&$\ddagger$&$\dagger$&$\ddagger$\\
			\hline
2&	\bf{0.22578}&	10&	0.60&	10&	0.64&	10&	0.67&	10&	3.38\\
3&	\bf{0.30676}&	10&	1.97&	10&	2.09&	10&	2.23&	10&	7.85\\
4&	\bf{0.37251}&	10&	4.38&	10&	4.84&	10&	5.17&	10&	13.63\\
5&	\bf{0.43555}&	10&	8.05&	10&	9.58&	10&	9.93&	10&	20.92\\
6&	\bf{0.48991}&	10&	13.89&	8&	15.60&	9&	16.15&	4&	28.96\\
7&	\bf{0.54401}&	10&	19.65&	10&	21.27&	10&	21.91&	9&	38.21\\
8&	\bf{0.59197}&	9&	27.82&	8&	31.59&	10&	34.34&	8&	49.26\\
9&	\bf{0.63683}&	10&	40.28&	10&	46.26&	10&	46.70&	9&	60.31\\
10&	\bf{0.67897}&	10&	52.37&	10&	55.28&	10&	59.14&	10&	73.48\\
11&	\bf{0.71825}&	10&	66.43&	10&	75.42&	9&	75.94&	9&	87.35\\
12&	\bf{0.75464}&	10&	80.79&	10&	91.16&	10&	96.64&	10&	100.59\\
13&	\bf{0.78744}&	10&	99.78&	10&	109.03&	9&	108.89&	10&	113.50\\
14&	\bf{0.81377}&	3&	117.91&	5&	148.88&	4&	136.75&	2&	125.12\\
15&	\bf{0.84165}&	7&	153.45&	8&	189.01&	7&	201.09&	4&	136.95\\
16&	\bf{0.86468}&	10&	176.79&	10&	189.31&	9&	188.24&	7&	145.95\\
17&	\bf{0.88340}&	9&	221.08&	8&	251.86&	6&	267.77&	4&	169.86\\
18&	\bf{0.90121}&	4&	231.24&	9&	239.83&	6&	294.76&	6&	195.90\\
19&	\bf{0.91541}&	6&	250.30&	7&	308.55&	2&	328.84&	5&	222.73\\
20&	\bf{0.92791}&	0&	297.46&	0&	378.29&	0&	449.72&	1&	251.16\\
21&	\bf{0.94052}&	6&	311.91&	4&	331.91&	2&	405.54&	3&	280.05\\
22&	\bf{0.95303}&	3&	356.73&	6&	423.16&	4&	438.65&	4&	312.48\\
\hline									
\multicolumn{2}{|l||}{Total}&		167&	2532.86&	173&	2923.56&	157&	3189.05&	145&	2437.64\\
\hline									
\multicolumn{2}{|l||}{Average:}&		7.95&	120.61&	8.24&	139.22&	7.48&	151.86&	6.90&	116.08\\
\hline									
\multicolumn{2}{|l||}{Percent$^*$}&\multicolumn{2}{|c||}{0.0313\%}	&\multicolumn{2}{|c||}{0.0251\%}&\multicolumn{2}{|c||}{0.0423\%}&\multicolumn{2}{|c|}{0.0857\%}\\
\hline
\multicolumn{10}{l}{$\dagger$ The number of times that the best known solution was found.}\\
\multicolumn{10}{l}{$\ddagger$ Total time in minutes for all ten runs.}\\
\multicolumn{10}{l}{$^*$ Average of 210 runs (10 runs each: $2\le p\le 22$) below the best known solution.}
		\end{tabular} 															
	\end{center}
\end{table}

\begin{table}[ht!]
	\begin{center}															
		\caption{\label{comp2}Comparing the Genetic Algorithm with Simulated Annealing ($p\ge 23$)}
		\medskip															\begin{tabular}{|c|c||c|c||c|c||c|c||c|c|}															\hline
			&Best&\multicolumn{2}{|c||}{Gen ($P=1$)}&\multicolumn{2}{|c||}{Gen ($P=2$)}&\multicolumn{2}{|c||}{Gen ($P=3$)}&\multicolumn{2}{|c|}{S.A.}\\
			\cline{3-10}
			$p$&Known&$\dagger$&$\ddagger$&$\dagger$&$\ddagger$&$\dagger$&$\ddagger$&$\dagger$&$\ddagger$\\
			\hline
23&	\bf{0.96196}&	0&	386.37&	2&	423.57&	0&	559.93&	3&	345.42\\
24&	\bf{0.96925}&	3&	415.35&	7&	466.89&	5&	500.73&	4&	380.03\\
25&	\bf{0.97566}&	6&	434.33&	3&	497.53&	5&	514.06&	3&	414.29\\
26&	\bf{0.98078}&	7&	433.50&	7&	533.18&	6&	570.72&	3&	447.90\\
27&	\bf{0.98564}&	2&	465.74&	5&	571.04&	2&	592.73&	3&	487.40\\
28&	\bf{0.98852}&	4&	493.45&	5&	675.23&	2&	694.79&	2&	524.36\\
29&	\bf{0.99085}&	1&	606.85&	1&	812.38&	1&	732.71&	3&	566.03\\
30&	\bf{0.99261}&	5&	657.98&	2&	735.25&	1&	805.24&	5&	604.02\\
31&	\bf{0.99430}&	4&	641.96&	4&	740.49&	2&	746.25&	1&	647.49\\
32&	\bf{0.99537}&	2&	666.92&	3&	705.93&	1&	797.02&	2&	692.48\\
33&	\bf{0.99633}&	0&	717.12&	0&	708.82&	1&	851.62&	0&	738.90\\
34&	\bf{0.99721}&	0&	645.29&	0&	829.53&	0&	888.11&	0&	788.46\\
35&	\bf{0.99803}$^*$&	0&	736.39&	0&	964.07&	0&	1020.54&	0&	841.77\\
36&	\bf{0.99859}&	0&	895.57&	0&	927.78&	0&	1391.82&	0&	892.83\\
37&	\bf{0.99903}&	0&	809.22&	0&	963.26&	0&	1521.67&	1&	945.97\\
38&	\bf{0.99945}$^*$&	0&	985.48&	2&	1063.38&	1&	1456.93&	0&	1000.55\\
39&	\bf{0.99969}&	0&	890.23&	0&	1023.01&	0&	1502.61&	0&	1054.54\\
40&	\bf{0.99987}&	0&	865.48&	0&	969.62&	0&	1333.65&	1&	1109.00\\
41&	\bf{0.99993}&	2&	827.42&	3&	930.79&	0&	1301.92&	0&	1163.53\\
42&	\bf{0.99997}&	0&	834.28&	1&	903.33&	1&	1235.90&	1&	1221.55\\
43&	\bf{1.00000}&	0&	817.42&	1&	922.30&	0&	1360.36&	0&	1281.10\\
44&	\bf{1.00000}&	9&	867.73&	10&	1019.97&	9&	1500.68&	0&	1344.23\\
\hline									
\multicolumn{2}{|l||}{Total}&		45&	15094.08&	56&	17387.35&	37&	21880.00&	32&	17491.84\\
\hline									
\multicolumn{2}{|l||}{Grand Total}&		212&	17626.94&	229&	20310.91&	194&	25069.06&	177&	19929.48\\
\hline									
\multicolumn{2}{|l||}{Grand Average:}&		4.93&	409.93&	5.33&	472.35&	4.51&	583.00&	4.12&	463.48\\
\hline
\multicolumn{2}{|l||}{\% of 220 runs$^+$}&\multicolumn{2}{|c||}{0.0251\%}	&\multicolumn{2}{|c||}{0.0228\%}&\multicolumn{2}{|c||}{0.0218\%}&\multicolumn{2}{|c|}{0.0264\%}\\
\hline
\multicolumn{2}{|l||}{\% of 430 runs$^+$}&\multicolumn{2}{|c||}{0.0276\%}	&\multicolumn{2}{|c||}{0.0238\%}&\multicolumn{2}{|c||}{0.0300\%}&\multicolumn{2}{|c|}{0.0501\%}\\
\hline
			\multicolumn{10}{l}{$\dagger$ The number of times that the best known solution was found.}\\
			\multicolumn{10}{l}{$\ddagger$ Total time in minutes for all ten runs.}\\
			\multicolumn{10}{l}{$^*$ A new best known solution.}\\
			\multicolumn{10}{l}{$^+$ \% of average  below the best known solution (10 runs each: $23\le p\le 44$;  $2\le p\le 44$).}
\end{tabular} 															
\end{center}
\end{table}

\subsection{Testing the Genetic Algorithm}

The best performing heuristic reported in \cite{DDK17} was simulated annealing \citep{KGV83}. By comparing the results with simulated annealing we implicitly also compare it with Tabu search and the ascent algorithm. We compare the proposed genetic algorithm to these results. We tested the genetic algorithm for $P=1,~2,~3$ second parents in the parents selection rule and $g=10,000$ iterations so that run times are similar to those reported in \cite{DDK17}.
We divided the 43 values of $p$ into two similar sized groups: $2\le p\le 22$ (21 problems) and $23\le p\le 44$ (22 problems).

In Table \ref{comp1} we report for the smaller values of $p$ the results for $P=1,2,3$ parents in the parents selection rule as well as the results by simulated annealing. For each method we also report: the number of times that the best known solution was found; the total time in minutes for all ten runs; and the percent of the average of 210 runs below the average of best known solutions. All three variants of the genetic algorithm performed better than the simulated annealing. Using $P=2$ parents in the parents selection rule provided the best results but run time was a bit longer than the time using $P=1$.

In Table \ref{comp2} we report the same results for the larger values of $p$ and also the grand total and average for all 43 problems. These results also indicate that the $P=2$ variant performed best.

\begin{figure}[ht!]
	\begin{center}
		\setlength{\unitlength}{1in}
		
		\begin{picture}(4,4.5)
		
		\includegraphics[width=4in]{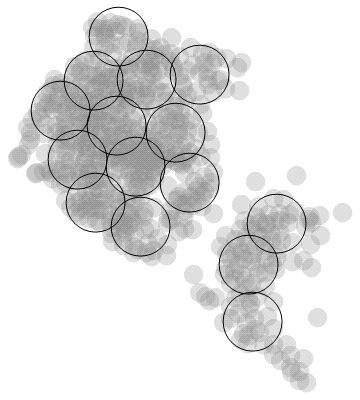}
		\put(-4,0){\bf Small shade: 1 mile radius around census tract}
		\put(-4,-0.2){\bf Large circle: 3 miles radius around towers}
		
		\end{picture}
		\vspace{0.1in}
		\caption{\label{15}Covering 85\% of Customers in Orange County by 15 Towers}
		
	\end{center}
\end{figure}

\subsection{Testing Location Anywhere in the Plane}

\begin{table}[htp!]
	\begin{center}															
		\caption{\label{comp100}Continuous solutions with 100 starting solutions}
		\medskip
		\small
		\begin{tabular}{|c|c||c|c||c|c|c|c||c|c|c|c|}
	\hline
	&&\multicolumn{2}{|c||}{SNOPT}&\multicolumn{4}{|c||}{N-M Random Starting Solutions}&\multicolumn{4}{|c|}{N-M Discrete Starting Solution}\\
	\cline{5-12}
	&Best&\multicolumn{2}{|c||}{~}&\multicolumn{2}{|c|}{$K=3$}&\multicolumn{2}{|c||}{$K=4$}&\multicolumn{2}{|c|}{$K=3$}&\multicolumn{2}{|c|}{$K=4$}\\
	\cline{3-12}
	$p$&Known&$\dagger$&$\ddagger$&$\dagger$&$\ddagger$&$\dagger$&$\ddagger$&$\dagger$&$\ddagger$&$\dagger$&$\ddagger$\\
	\hline
2&	{\bf 0.22680}&	6.13\%&	0.68&		0.00\%&	0.04&			0.00\%&	0.05&							0.00\%&	0.03&							0.00\%&	0.04\\
3&	{\bf 0.30888}&	13.97\%&	1.12&		0.01\%&	0.12&			0.02\%&	0.17&							0.00\%&	0.08&							0.00\%&	0.10\\
4&	{\bf 0.37597}&	18.68\%&	1.83&		0.01\%&	0.26&			0.01\%&	0.34&							0.00\%&	0.17&							0.00\%&	0.23\\
5&	{\bf 0.44094}&	15.99\%&	2.44&		0.01\%&	0.40&			0.00\%&	0.63&							0.00\%&	0.28&							0.00\%&	0.34\\
6&	{\bf 0.49555}&	15.37\%&	2.94&		0.02\%&	0.67&			0.00\%&	0.88&							0.28\%&	0.50&							0.29\%&	0.68\\
7&	{\bf 0.54924}&	14.49\%&	4.19&		0.08\%&	1.05&			0.02\%&	1.25&							0.07\%&	0.69&							0.08\%&	0.87\\
8&	{\bf 0.59883}&	19.62\%&	4.69&		0.03\%&	1.37&			0.04\%&	1.81&							0.00\%&	0.90&							0.00\%&	1.15\\
9&	{\bf 0.64359}&	17.11\%&	5.01&		0.06\%&	1.92&			0.06\%&	2.56&							0.02\%&	1.28&							0.01\%&	1.87\\
10&	{\bf 0.68718}&	18.06\%&	5.92&		0.04\%&	2.50&			0.06\%&	3.48&							0.00\%&	1.67&							0.00\%&	2.33\\
11&	{\bf 0.72856}&	21.74\%&	7.30&		0.02\%&	3.18&			0.01\%&	4.19&							0.00\%&	1.96&							0.00\%&	2.62\\
12&	{\bf 0.76548}&	21.02\%&	7.43&		0.17\%&	3.84&			0.18\%&	4.94&							0.00\%&	2.30&							0.00\%&	3.09\\
13&	{\bf 0.79850}&	18.35\%&	10.32&		0.03\%&	4.28&			0.10\%&	5.74&							0.00\%&	2.76&							0.00\%&	3.70\\
14&	{\bf 0.82320}&	19.06\%&	10.05&		0.07\%&	5.22&			0.06\%&	7.24&							0.05\%&	2.82&							0.05\%&	3.59\\
15&	{\bf 0.84941}&	14.98\%&	10.74&		0.26\%&	5.66&			0.12\%&	7.71&							0.00\%&	2.86&							0.01\%&	3.91\\
16&	{\bf 0.87492}&	22.27\%&	12.01&		0.22\%&	6.53&			0.02\%&	8.70&							0.00\%&	3.64&							0.00\%&	5.23\\
17&	{\bf 0.89310}&	21.96\%&	13.43&		0.12\%&	7.07&			0.06\%&	9.79&							0.02\%&	3.67&							0.02\%&	5.32\\
18&	{\bf 0.91138}&	17.74\%&	13.44&		0.04\%&	8.16&			0.23\%&	10.70&							0.02\%&	4.01&							0.01\%&	5.56\\
19&	{\bf 0.92767}&	20.02\%&	14.85&		0.52\%&	8.80&			0.53\%&	11.70&							0.01\%&	4.48&							0.01\%&	5.78\\
20&	{\bf 0.93940}&	17.45\%&	14.48&		0.23\%&	9.66&			0.76\%&	12.33&							0.02\%&	4.80&							0.01\%&	6.55\\
21&	{\bf 0.95053}&	15.75\%&	17.09&		0.35\%&	9.20&			0.10\%&	13.45&							0.17\%&	5.97&							0.17\%&	8.92\\
22&	{\bf 0.96060}&	16.68\%&	17.23&		0.06\%&	10.71&			0.20\%&	14.65&							0.00\%&	6.35&							0.00\%&	8.83\\
23&	{\bf 0.96917}&	14.86\%&	18.47&		0.23\%&	10.41&			0.10\%&	14.65&							0.03\%&	7.60&							0.03\%&	10.94\\
24&	{\bf 0.97711}&	14.06\%&	19.11&		0.71\%&	10.42&			0.29\%&	15.24&							0.01\%&	6.66&							0.02\%&	9.30\\
25&	{\bf 0.98309}&	15.69\%&	19.74&		0.19\%&	11.87&			0.14\%&	16.10&							0.03\%&	7.78&							0.02\%&	10.34\\
26&	{\bf 0.98770}&	13.78\%&	22.55&		0.22\%&	11.57&			0.21\%&	15.14&							0.04\%&	8.33&							0.04\%&	10.81\\
27&	{\bf 0.99210}&	14.98\%&	22.84&		0.36\%&	11.36&			0.37\%&	16.41&							0.02\%&	8.26&							0.02\%&	11.90\\
28&	{\bf 0.99377}&	8.67\%&	25.02&		0.08\%&	11.45&			0.19\%&	16.77&							0.01\%&	6.82&							0.00\%&	9.08\\
29&	{\bf 0.99604}&	12.29\%&	25.82&		0.16\%&	11.70&			0.21\%&	16.22&							0.01\%&	7.02&							0.01\%&	9.99\\
30&	{\bf 0.99704}&	13.53\%&	23.36&		0.12\%&	11.48&			0.01\%&	15.17&							0.00\%&	7.12&							0.00\%&	10.34\\
31&	{\bf 0.99811}&	10.05\%&	25.37&		0.03\%&	10.87&			0.09\%&	16.43&							0.00\%&	5.83&							0.00\%&	8.60\\
32&	{\bf 0.99857}&	12.11\%&	28.51&		0.08\%&	11.40&			0.11\%&	15.67&							0.01\%&	5.51&							0.00\%&	7.51\\
33&	{\bf 0.99947}&	10.93\%&	26.53&		0.08\%&	10.84&			0.10\%&	14.99&							0.02\%&	5.37&							0.02\%&	7.71\\
34&	{\bf 0.99963}&	8.64\%&	28.57&		0.07\%&	9.59&			0.13\%&	14.07&							0.02\%&	4.45&							0.02\%&	6.11\\
35&	{\bf 0.99998}&	10.14\%&	29.58&		0.09\%&	10.91&			0.10\%&	13.14&							0.00\%&	4.78&							0.00\%&	6.82\\
36&	{\bf 1.00000}&	9.44\%&	29.33&		0.04\%&	9.24&			0.06\%&	13.45&							0.00\%&	3.65&							0.00\%&	5.47\\
37&	{\bf 1.00000}&	7.26\%&	30.80&		0.01\%&	10.17&			0.02\%&	12.29&							0.00\%&	2.54&							0.00\%&	3.32\\
38&	{\bf 1.00000}&	7.55\%&	31.97&		0.10\%&	8.91&			0.00\%&	12.33&							0.00\%&	1.96&							0.00\%&	2.69\\
39&	{\bf 1.00000}&	8.48\%&	32.32&		0.03\%&	9.71&			0.01\%&	13.75&							0.00\%&	1.45&							0.00\%&	1.85\\
40&	{\bf 1.00000}&	6.09\%&	33.75&		0.00\%&	9.30&			0.01\%&	11.88&							0.00\%&	1.31&							0.00\%&	1.91\\
\hline																										
\multicolumn{2}{|l||}{Average:}&		14.49\%&	16.69&		0.126\%&	7.23&			0.121\%&	9.90&							0.022\%&	3.79&							0.022\%&	5.27\\
			\hline
			\multicolumn{12}{l}{$\dagger$ Percent cover below best known solution.}\\
			\multicolumn{12}{l}{$\ddagger$ Total time in minutes for all runs.}\\
		\end{tabular} 															
	\end{center}
\end{table}

\begin{table}[htp!]
	\begin{center}															
		\caption{\label{comp1000}Continuous solutions with 1000 starting solutions}
		\medskip
		\small
		\setlength{\tabcolsep}{5pt}
		\begin{tabular}{|c|c||c|c||c|c|c|c||c|c|c|c|}
\hline
&&\multicolumn{2}{|c||}{SNOPT}&\multicolumn{4}{|c||}{N-M Random Starting Solutions}&\multicolumn{4}{|c|}{N-M Discrete Starting Solution}\\
\cline{5-12}
&Best&\multicolumn{2}{|c||}{~}&\multicolumn{2}{|c|}{$K=3$}&\multicolumn{2}{|c||}{$K=4$}&\multicolumn{2}{|c|}{$K=3$}&\multicolumn{2}{|c|}{$K=4$}\\
			\cline{3-12}
			$p$&Known&$\dagger$&$\ddagger$&$\dagger$&$\ddagger$&$\dagger$&$\ddagger$&$\dagger$&$\ddagger$&$\dagger$&$\ddagger$\\
			\hline
2&	{\bf 0.22680}&	2.09\%&	7.16&		0.00\%&	0.42&			0.00\%&	0.56&							0.00\%&	0.29&							0.00\%&	0.38\\
3&	{\bf 0.30888}&	5.42\%&	11.38&		0.00\%&	1.24&			0.00\%&	1.65&							0.00\%&	0.83&							0.00\%&	1.06\\
4&	{\bf 0.37597}&	11.41\%&	17.25&		0.00\%&	2.62&			0.00\%&	3.47&							0.00\%&	1.75&							0.00\%&	2.27\\
5&	{\bf 0.44094}&	10.24\%&	24.29&		0.00\%&	4.49&			0.00\%&	5.85&							0.00\%&	2.69&							0.00\%&	3.49\\
6&	{\bf 0.49555}&	12.67\%&	31.18&		0.00\%&	7.14&			0.00\%&	9.31&							0.23\%&	4.96&							0.26\%&	6.58\\
7&	{\bf 0.54924}&	12.98\%&	38.51&		0.00\%&	10.10&			0.02\%&	13.19&							0.07\%&	6.63&							0.02\%&	9.03\\
8&	{\bf 0.59883}&	15.65\%&	46.07&		0.01\%&	13.89&			0.01\%&	18.82&							0.00\%&	9.09&							0.00\%&	12.02\\
9&	{\bf 0.64359}&	13.98\%&	54.89&		0.05\%&	18.55&			0.05\%&	25.43&							0.00\%&	13.35&							0.01\%&	18.46\\
10&	{\bf 0.68718}&	15.49\%&	63.75&		0.02\%&	24.58&			0.01\%&	33.37&							0.00\%&	16.82&							0.00\%&	22.39\\
11&	{\bf 0.72856}&	17.13\%&	76.92&		0.00\%&	30.66&			0.01\%&	41.60&							0.00\%&	19.78&							0.00\%&	26.65\\
12&	{\bf 0.76548}&	15.26\%&	85.26&		0.01\%&	37.64&			0.01\%&	49.27&							0.00\%&	22.55&							0.00\%&	30.54\\
13&	{\bf 0.79850}&	15.79\%&	89.34&		0.03\%&	43.01&			0.01\%&	57.50&							0.00\%&	27.28&							0.00\%&	37.46\\
14&	{\bf 0.82320}&	15.52\%&	103.02&		0.00\%&	51.20&			0.01\%&	69.08&							0.05\%&	28.15&							0.05\%&	36.79\\
15&	{\bf 0.84941}&	18.88\%&	109.40&		0.00\%&	56.78&			0.00\%&	76.83&							0.00\%&	29.18&							0.00\%&	39.17\\
16&	{\bf 0.87492}&	16.02\%&	123.56&		0.00\%&	64.34&			0.02\%&	87.75&							0.00\%&	37.02&							0.00\%&	51.28\\
17&	{\bf 0.89310}&	14.53\%&	135.86&		0.02\%&	72.16&			0.00\%&	99.00&							0.02\%&	37.75&							0.02\%&	52.11\\
18&	{\bf 0.91138}&	18.53\%&	139.27&		0.00\%&	80.44&			0.02\%&	109.25&							0.01\%&	39.90&							0.01\%&	54.68\\
19&	{\bf 0.92767}&	12.90\%&	151.05&		0.25\%&	85.98&			0.07\%&	119.69&							0.00\%&	43.45&							0.00\%&	58.93\\
20&	{\bf 0.93940}&	15.99\%&	166.43&		0.00\%&	95.51&			0.07\%&	128.35&							0.01\%&	47.59&							0.01\%&	66.72\\
21&	{\bf 0.95053}&	15.96\%&	173.13&		0.19\%&	98.13&			0.00\%&	136.20&							0.17\%&	60.85&							0.17\%&	84.66\\
22&	{\bf 0.96060}&	14.57\%&	178.20&		0.06\%&	107.45&			0.09\%&	143.36&							0.00\%&	65.55&							0.00\%&	92.09\\
23&	{\bf 0.96917}&	10.26\%&	192.27&		0.00\%&	108.61&			0.02\%&	151.42&							0.02\%&	74.32&							0.02\%&	103.48\\
24&	{\bf 0.97711}&	13.20\%&	202.90&		0.00\%&	109.07&			0.03\%&	152.71&							0.00\%&	67.48&							0.00\%&	94.54\\
25&	{\bf 0.98309}&	12.39\%&	209.77&		0.00\%&	112.80&			0.14\%&	158.04&							0.02\%&	76.83&							0.02\%&	107.45\\
26&	{\bf 0.98770}&	11.80\%&	221.78&		0.00\%&	115.21&			0.03\%&	157.99&							0.02\%&	79.90&							0.03\%&	113.28\\
27&	{\bf 0.99210}&	8.54\%&	230.00&		0.02\%&	115.02&			0.22\%&	164.06&							0.00\%&	84.58&							0.02\%&	115.98\\
28&	{\bf 0.99377}&	9.94\%&	236.59&		0.06\%&	118.09&			0.04\%&	165.25&							0.00\%&	64.14&							0.00\%&	90.03\\
29&	{\bf 0.99604}&	11.47\%&	254.07&		0.08\%&	114.07&			0.12\%&	163.13&							0.00\%&	70.34&							0.01\%&	100.62\\
30&	{\bf 0.99704}&	8.96\%&	260.81&		0.04\%&	113.23&			0.01\%&	162.62&							0.00\%&	70.05&							0.00\%&	99.71\\
31&	{\bf 0.99811}&	10.80\%&	266.64&		0.00\%&	110.50&			0.08\%&	156.56&							0.00\%&	58.93&							0.00\%&	83.28\\
32&	{\bf 0.99857}&	5.86\%&	274.88&		0.07\%&	110.27&			0.04\%&	152.58&							0.00\%&	53.28&							0.00\%&	77.36\\
33&	{\bf 0.99947}&	7.09\%&	295.22&		0.02\%&	102.86&			0.07\%&	147.00&							0.00\%&	54.14&							0.02\%&	75.36\\
34&	{\bf 0.99963}&	7.03\%&	294.63&		0.03\%&	100.49&			0.01\%&	139.83&							0.00\%&	48.19&							0.02\%&	65.98\\
35&	{\bf 0.99998}&	6.28\%&	305.51&		0.01\%&	100.78&			0.05\%&	134.79&							0.00\%&	46.51&							0.00\%&	67.91\\
36&	{\bf 1.00000}&	7.11\%&	312.73&		0.02\%&	94.91&			0.01\%&	136.91&							0.00\%&	37.58&							0.00\%&	53.99\\
37&	{\bf 1.00000}&	5.76\%&	325.34&		0.01\%&	93.50&			0.01\%&	132.18&							0.00\%&	24.67&							0.00\%&	33.30\\
38&	{\bf 1.00000}&	5.44\%&	337.85&		0.00\%&	92.82&			0.00\%&	124.58&							0.00\%&	19.85&							0.00\%&	27.05\\
39&	{\bf 1.00000}&	4.68\%&	339.68&		0.00\%&	91.68&			0.00\%&	124.43&							0.00\%&	14.52&							0.00\%&	19.33\\
40&	{\bf 1.00000}&	5.35\%&	350.84&		0.00\%&	88.25&			0.00\%&	122.49&							0.00\%&	13.73&							0.00\%&	18.80\\
\hline																										
\multicolumn{2}{|l||}{Average:}&		11.36\%&	172.75&		0.026\%&	71.76&			0.033\%&	99.39&							0.016\%&	37.81&							0.018\%&	52.67\\
\hline																									\multicolumn{12}{l}{$\dagger$ Percent cover below best known solution.}\\
			\multicolumn{12}{l}{$\ddagger$ Total time in minutes for all runs.}\\
		\end{tabular} 															
	\end{center}
\end{table}

\begin{figure}[ht!]
	\begin{center}
		\setlength{\unitlength}{1in}
\hspace{-3in}
		\begin{picture}(3,2.8)
		\includegraphics[width=3in]{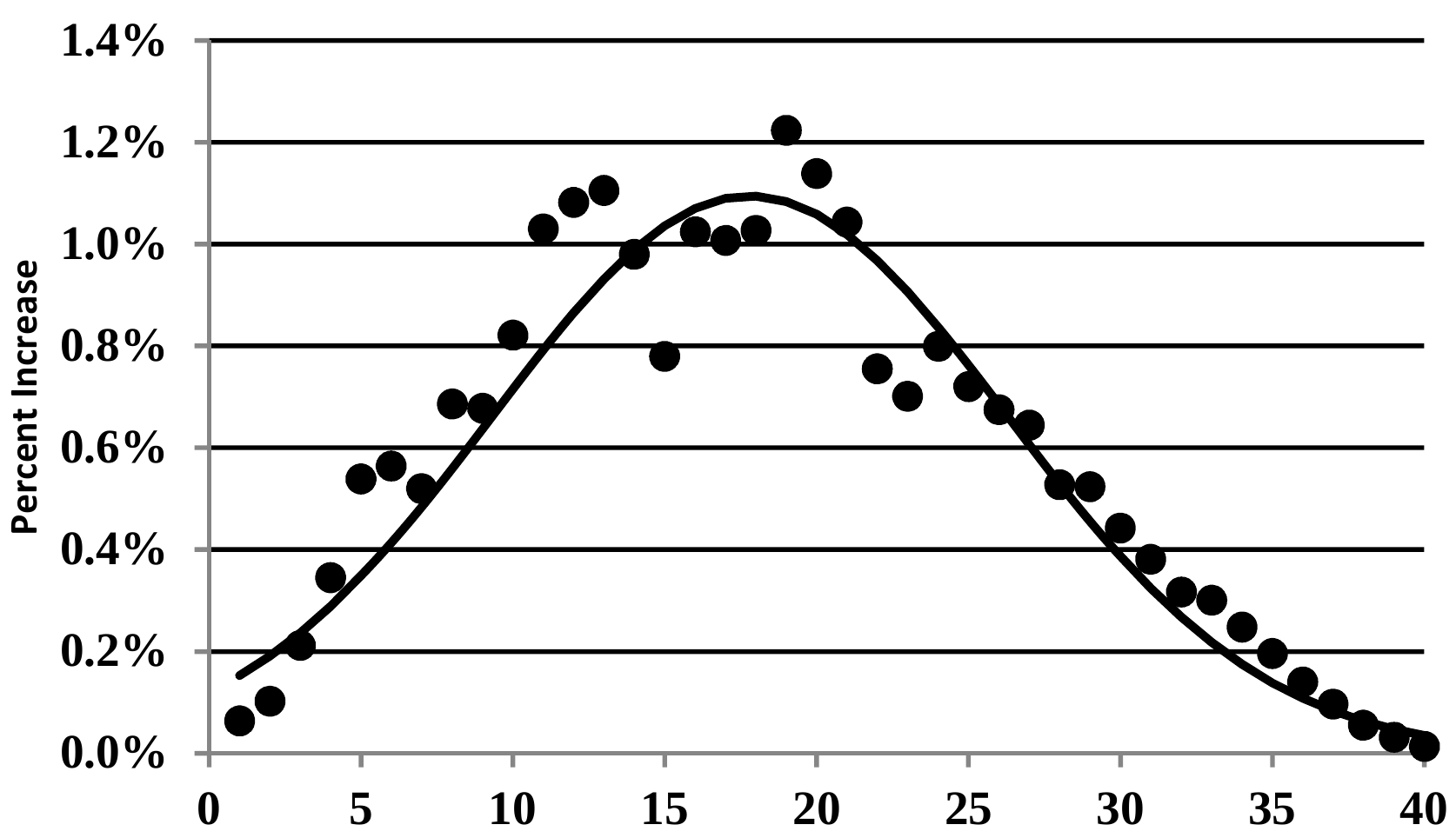}
			\includegraphics[width=3in]{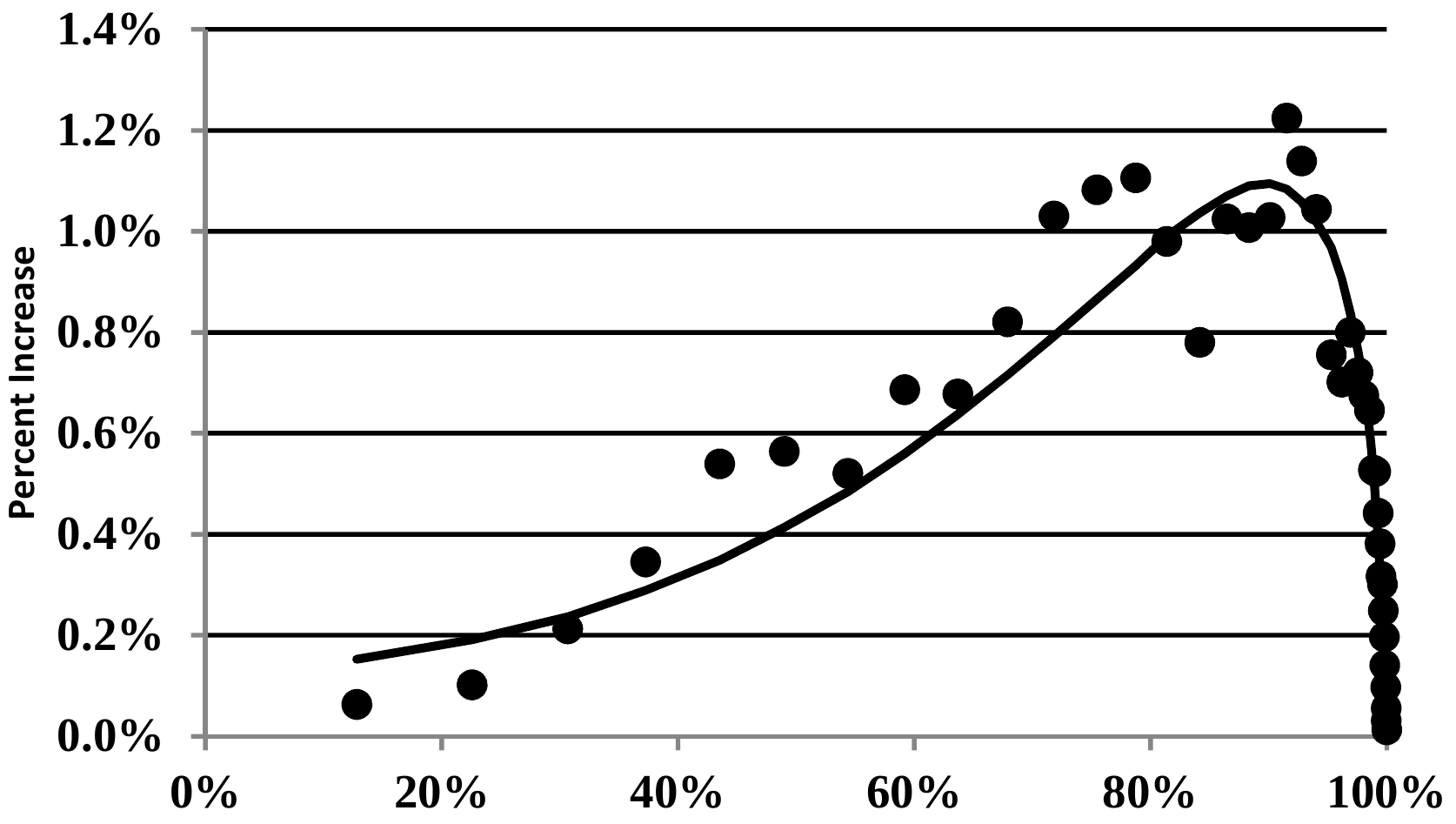}
\put(-5,-0.15){\bf\scriptsize Number of Towers}
\put(-2.,-0.15){\bf\scriptsize Discrete Percent Cover}
		\end{picture}
		\bigskip
		\caption{\label{fig7}Cover Increase at the Continuous Location}
	\end{center}
\end{figure}

We compared maximum cover of all $n=577$ census tracts. SNOPT (using the hexagonal pattern) and Nelder-Mead (using Gaussian quadrature) )were run from random starting solutions. In order to improve some best known solutions we also report results by Nelder-Mead starting from the discrete best known solution found by the genetic algorithm and in \cite{DDK17}. Each variant was run for 100 and 1000 replications for each $2\le p\le 40$.

The results are summarized in Tables \ref{comp100} and \ref{comp1000}. Run times by Nelder-Mead are quite short. The average run time for one run is about 5 seconds from a random starting solution and 3 seconds from the best known discrete starting solution (it was repeated 100 or 1000 times). SNOPT was slower and resulted in much poorer results and thus it is not discussed further. The random starting solution performed better for small values of $p$ and the discrete starting solution performed better for larger values of $p$. The $K=3$ yielded comparable results to those obtained for $K=4$ in a slightly shorter run time.  Full cover was obtained for $p\ge 36$ instances. Cover of 99.999\% which is practically full cover was obtained for $p=35$. In Figure~\ref{15} the best cover obtained by 15 towers is depicted. Full cover will be obtained if all gray shaded areas are covered by towers within 3 miles.

In Figure \ref{fig7} we compare the best known cover for the discrete problem with the best known cover for location anywhere in the plane.
As would be expected, the cover by locating anywhere in the plane is higher than the discrete cover. In some cases more than 1\% is added to the cover (out of 100\%).   On the left panel we show the increase as a function of the the number of towers.  We added a best fit function. On the right panel we show the same values (and the best fit graph) with the discrete cover on the $x$-axis. Additional cover is low for few facilities because the facilities in the continuous solution are located quite close to the demand points, and it is low for a large number of facilities because there is not much extra market share remaining. For our case study the maximum additional market share is obtained around 90\% discrete cover.

\begin{figure}[ht!]
	\begin{center}
		\setlength{\unitlength}{2in}
		\caption{\label{NOrange}Covering All Customers in North Orange County}
		\begin{picture}(3,1.4)
		
		\includegraphics[width=6in]{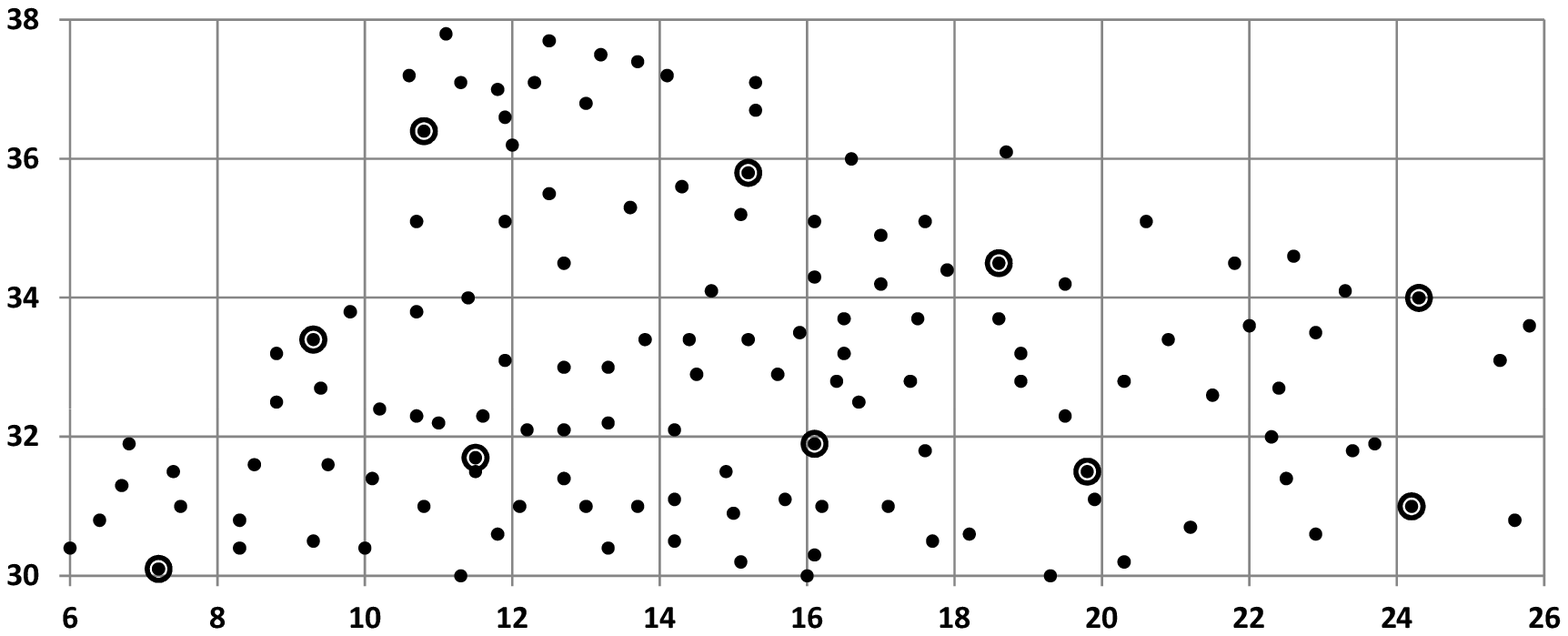}
		\put(-2.22,1.22){\bf Ten Towers Located at Census Tracts}
		\put(-3.05,0.7){$y$}
		\put(-1.6,-0.05){$x$}
		\end{picture}
		\begin{picture}(3,1.4)
		\thicklines
		\includegraphics[width=6in]{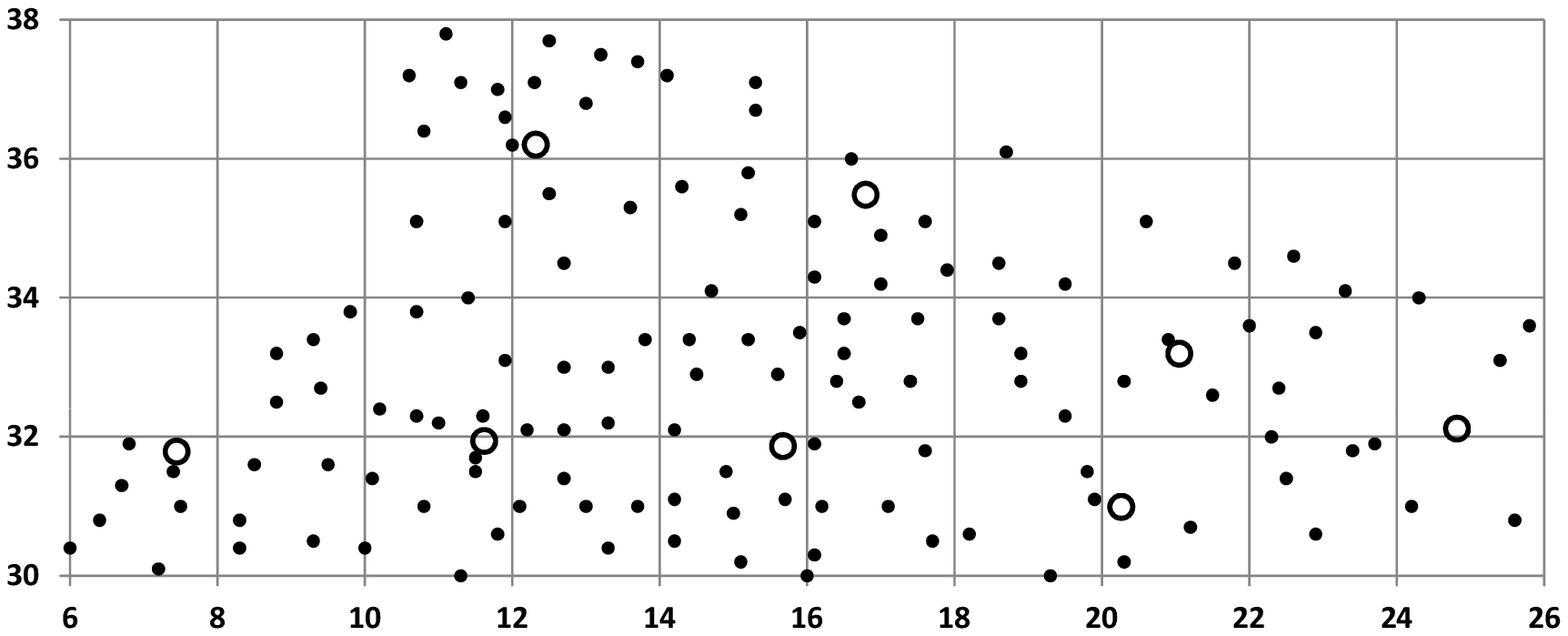}
		\put(-2.3,1.22){\bf Eight Towers Located Anywhere in the Area}
		\put(-2.05,-0.1){\circle*{0.03}}
		\put(-2.0,-0.125){\bf Census tract;}
		\put(-1.4,-0.1){\circle{0.04}}
		\put(-1.4,-0.1){\circle{0.035}}
		\put(-1.35,-0.125){\bf Tower.}
		\put(-3.05,0.7){$y$}
		\put(-1.6,-0.05){$x$}
		
		\end{picture}
	\end{center}
\end{figure}

We also investigated covering North Orange County  ($n=13$1 census tracts for $y\ge 30$). Full cover for North Orange County was found by Nelder-Mead for $p\ge 8$ facilities. SNOPT was not tested on North Orange County in view of its performance for all Orange County. Nelder Mead was run 100 times and full cover obtained 18 times for $p=8$. Running time for all 100 runs for $p=8$ was just 5.8 seconds. When the locations are restricted to census tracts full cover requires $p=10$ facilities \citep{DDK17}. In earlier gradual cover models \citep{DrDr14,BDK11b} 13 towers are required for full cover. The solutions are depicted in Figure \ref{NOrange}. Clearly, locating anywhere in the plane  provides much more flexibility.

\section{\label{concl}Conclusions}

We investigated locating a given number of facilities providing maximum cover of a set of demand points applying the directional gradual cover model. This is a realistic and practical model compared with other gradual cover formulations. It is, however, more difficult to solve. A better heuristic approach is obtained for the discrete case by a specially designed genetic algorithm.

We also find the best locations anywhere in the plane by Nelder-Mead and SNOPT.
Nelder-Mead provided superior results because the objective function has discontinuous derivatives. All purpose non-linear optimization solvers perform poorly because of this discontinuity.

Some interesting and useful extensions can be investigated in future research. For example,
\begin{itemize}
\item 
 The conditional version of the problem \citep{Mi80,BS90,OZ02}. This means that several facilities already exist in the area and additional facilities need to be located to maximize total cover.
 \item The problem can be defined on the globe with a large cover radius determined by a limited flight time. The distance required for a given flight time can be in a range due to variable flight conditions.
 \item The cost of a facility is a known function of the cover radius. For example, the cost consists of a set-up cost plus a variable cost which is proportional to the square of the radius. A budget is available for constructing new facilities. The cover radii of facilities are additional variables in the model. The conditional version of this problem is also of interest.
\end{itemize}

{\centering \section*{Appendix: Finding the Joint Cover of One Demand Point by Hexagonal Pattern}}

The points in the hexagonal pattern are defined by two sequences (all the combinations of the two lists for $x$ and $y$): $(x=0,\pm1,\pm2,\ldots,y=0,\pm\sqrt3,\pm2\sqrt3,\ldots)$ and $(x=\pm\frac12,\pm\frac32,\pm\frac52,\ldots,y=\pm\frac{\sqrt3}2,\pm\frac{3\sqrt3}2\pm\frac{5\sqrt3}2,\ldots)$.
All the points with distances up to $\sqrt{M}$ for some $M$ are selected. $M=52$ yields $N=199$ points and $M=110, 220$ yield $N=397, 805$ points. The points are then adjusted by a factor of  $\sqrt{\frac{2\pi }{N\sqrt3}}$.

Every demand point is represented by $N$ hexagonal pattern points (see for example Figure \ref{hex}). The distance between each of the $N$ points to all facilities is calculated to find out whether the point is covered or not. The total number of points covered, $TC$, is found and the cover estimate is $\frac{TC}{N}$.

Time saving measures can be implemented. We first find the distances $d_j$ between the demand point (the center of the circle) and facility $j$ for $j=1,\ldots p$. If there is a facility for which $d_j\le D_j-R_i$, the whole circle is covered. The cover is ``1", and there is no need to check the other $N-1$ points. If  $d_j\ge D_j+R_i$, facility $j$ provides no cover and can be removed from the list of facilities to be checked by the other $N-1$ points. A list, possibly shorter than $p$ members, of facilities is created. Only facilities in this list need to be checked for each of the remaining $N-1$ points.

Note that once a facility covering a point is identified, the point is covered and there is no need to check additional facilities. A point is not covered if and only if all facilities  do not cover it.

\renewcommand{\baselinestretch}{1}
\renewcommand{\arraystretch}{1}
\large
\normalsize
\small

\bibliographystyle{apalike}

\bigskip


\end{document}